\documentclass[a4paper, 11pt]{article}

\usepackage{longtable}
\usepackage{geometry}
\usepackage{graphicx}
\usepackage{subfig}
\usepackage{color}
\usepackage{xcolor}
\usepackage{bm}
\usepackage{fixmath}
\usepackage{mathrsfs}
\usepackage{ulem}
\usepackage{multirow}
\usepackage{pdflscape}
\usepackage{rotating}
\usepackage{lscape}
\usepackage{pdflscape}
\usepackage{titlesec}
\usepackage{float}
\usepackage{hyperref}

\titleclass{\subsubsubsection}{straight}[\subsection]

\newcounter{subsubsubsection}[subsubsection]
\renewcommand\thesubsubsubsection{\thesubsubsection.\arabic{subsubsubsection}}

\titleformat{\subsubsubsection}
  {\normalfont\normalsize\bfseries}{\thesubsubsubsection}{1em}{}
\titlespacing*{\subsubsubsection}
{0pt}{3.25ex plus 1ex minus .2ex}{1.5ex plus .2ex}

\makeatletter

\renewcommand\paragraph{\@startsection{paragraph}{5}{\z@}%
  {3.25ex \@plus1ex \@minus.2ex}%
  {-1em}%
  {\normalfont\normalsize\bfseries}}

\renewcommand\subparagraph{\@startsection{subparagraph}{6}{\parindent}%
  {3.25ex \@plus1ex \@minus .2ex}%
  {-1em}%
  {\normalfont\normalsize\bfseries}}

\def\toclevel@subsubsubsection{4}
\def\toclevel@paragraph{5}
\def\toclevel@paragraph{6}
\def\l@subsubsubsection{\@dottedtocline{4}{7em}{4em}}
\def\l@paragraph{\@dottedtocline{5}{10em}{5em}}
\def\l@subparagraph{\@dottedtocline{6}{14em}{6em}}

\makeatother

\def\b{\begin{eqnarray}}
\def\e{\end{eqnarray}}

\def\n{\noindent}

\def\1{\vskip.1cm}
\def\2{\vskip.2cm}
\def\3{\vskip.3cm}
\def\5{\vskip.5cm}

\def\f{$f_n(x)\,$}
\def\g{$f_{n-1}(x)\,$}

\def\a{\alpha}

\makeatletter
\newcommand*\bigdot{\mathpalette\bigdot@{.5}}
\newcommand*\bigdot@[2]{\mathbin{\vcenter{\hbox{\scalebox{#2}{$\m@th#1\bullet$}}}}}
\makeatother

\begin{document}

\begin{center}
{\huge \textbf{New Bounds on the Real Polynomial \3 Roots}}

\vspace{9mm}
\noindent
{\large \bf Emil M. Prodanov} \vskip.4cm
{\it School of Mathematical Sciences, Technological University Dublin,
\vskip.1cm
City Campus, Kevin Street, Dublin, D08 NF82, Ireland,}
\vskip.1cm
{\it E-Mail: emil.prodanov@tudublin.ie} \\
\vskip.5cm
\end{center}

\vskip2cm

\n
\begin{abstract}
\n
The presented analysis determines several new bounds on the roots of the equation $a_n x^n + a_{n-1} x^{n-1} + \cdots + a_0 = 0$ (with $a_n > 0$). All proposed new bounds are lower than the Cauchy bound max$\{1, \sum_{j=0}^{n-1} |a_j/a_n| \}$. Firstly, the Cauchy bound formula is derived by presenting it in a new light --- through a recursion. It is shown that this recursion could be exited at earlier stages and, the earlier the recursion is terminated, the lower the resulting root bound will be. Following a separate analysis, it is further demonstrated that a significantly lower root bound can be found if the summation in the Cauchy bound formula is made not over each one of the coefficients $a_0, a_1, \ldots, a_{n-1}$, but only over the negative ones. The sharpest root bound in this line of analysis is shown to be the larger of 1 and the sum of the absolute values of all negative coefficients of the equation divided by the largest positive coefficient. The following bounds are also found in this paper: $\mbox{max}\{ 1, ( \sum_{j = 1}^{q} B_j/A_l )^{1/(l-k)}\}$, where $B_1, B_2, \ldots B_q$ are the absolute values of all of the negative  coefficients in the equation, $k$ is the highest degree of a monomial with a negative coefficient, $A_l$ is the positive coefficient of the term $A_l x^l$ for which $k< l \le n$.

\end{abstract}


\newpage

\n
The roots of the equation $a_n x^n + a_{n-1} x^{n-1} + \cdots + a_0 = 0$ are bound from above by the unique positive root of the associated Cauchy polynomial $|a_n| x^n - |a_{n-1}| x^{n-1} - \cdots - |a_0|$. The Cauchy formula yields that the upper bound of the roots of the equation is max$\{1, \sum_{j=0}^{n-1} |a_j/a_n| \}$ --- see section 8.1 in \cite{rs}. In the first part of the analysis presented in this paper, this formula is derived  from a different perspective --- through a recursion, --- following the idea of splitting polynomials into two parts and studying the ``interaction" between these parts \cite{1} (i.e. studying the intersection points of their graphs) --- a method successfully used also for the full classification of the roots of the cubic \cite{1} and the quartic equation \cite{2} in terms of the equation coefficients. This recursion involves bounding the unique positive root of a particular equation with the unique positive root of a subsidiary equation of degree one less. The recursion ends with the determination of the root of a linear equation and this root is exactly $\sum_{j=0}^{n-1} |a_j/a_n|$. If, instead, the recursion is terminated at an earlier stage --- that of a quadratic, cubic, or quartic equation ---  the resulting root bound will be lower, as shown in this work. All of these analytically determinable new bounds are lower than the Cauchy bound. \\
It is separately shown, following a different line of analysis, that one does not have to sum over all coefficients $a_j$ in $\sum_{j=0}^{n-1} |a_j/a_n|$, but only over the negative ones. This results in a significantly lower root bound that the Cauchy bound. It is demonstrated further that this new bound can be made even lower by finding a denominator, greater than $|a_n|$. The sharpest upper bound of the roots of the general polynomial $a_n x^n + a_{n-1} x^{n-1} + \cdots + a_0$ with $a_n > 0$ that can be found following this analysis, is either 1 or the smallest of the unique positive roots of all Cauchy polynomials that can be extracted form this polynomial with preservation of all terms with negative coefficients --- whichever is larger. The latter is the sum of the absolute values of all negative coefficients of the equation divided by the largest positive coefficient. \\
The following bounds are also found in this paper: $\mbox{max}\{ 1, ( \sum_{j = 1}^{q}  B_j/A_l  )^{1/(l-k) }\}$, where $B_1, B_2, \ldots B_q$ are the absolute values of all of the negative  coefficients in the equation, $k$ is the highest degree of a monomial with a negative coefficient, $A_l$ is the positive coefficient of the term $A_l x^l$ with $k< l \le n$. \\
The lower bound on the unique positive root of a Cauchy polynomial is also determined.
\vskip1cm
\n
A {\it Cauchy polynomial} of degree $n$ is defined \cite{rs} as a polynomial in $x$ such that the coefficient of $x^n$ is positive and the coefficients of all of its remaining terms --- negative. That is, the Cauchy polynomial has the form $|a_n| x^n - |a_{n-1}| x^{n-1} - \cdots - |a_0|$. As there must be at least one term with negative coefficient, a monomial cannot be a Cauchy polynomial. \\
The Cauchy polynomial has, by Descartes' rule of signs, a unique positive root, say $\mu$ (as there is only one sign change in the sequence of its coefficients), and, by Cauchy's theorem \cite{rs}, the root $\mu$ provides the Cauchy bound on the roots of the general polynomial equation of degree $n$
\b
\label{eqn}
a_n x^n + a_{n-1} x^{n-1} + \cdots + a_0 = 0.
\e
In this equation $a_n \ne 0$ and, without losing generality, it will be assumed further that $a_n > 0$. \\
If all of the coefficients of this equation have the same sign, then it will not have positive roots and therefore 0 will be the upper bound of the roots (if real roots exist). This special case will not be considered further. \\
Before addressing the general equation (\ref{eqn}), re-write the Cauchy polynomial as $\a_n x^n - \a_{n-1} x^{n-1} - \cdots - \a_0$, where all coefficients $\a_j$ are positive (the number of terms $\a_j x^j$ for which $j < n$ may be between 1 and $n$), and consider the associated equation
\b
\label{stump}
\a_n x^n - \a_{n-1} x^{n-1} - \cdots - \a_0 = 0.
\e
As the unique positive root of this equation is sought, one can assume, without loss of generality, that $\a_0 \ne 0$. (If $\a_0$ happens to be zero, one identifies 0 as a root, reduces everywhere the powers of $x$ and the indexes of the coefficients by one unit and arrives at an equation of the same type, but of degree $n-1$. If the ``new" $\a_0$ also happens to be zero, this procedure should be repeated until a non-zero coefficient is encountered --- it is guaranteed to exist by the definition.) \\
Equation (\ref{stump}) can be viewed as
\b
\label{f}
f_n(x) = \a_0,
\e
where $f_n(x) \equiv \a_n x^n - \a_{n-1} x^{n-1} - \cdots - \a_1 x$, and also as
\b
\label{g}
x f_{n-1}(x) = \a_0,
\e
where $f_{n-1}(x) \equiv \a_n x^{n-1} - \a_{n-1} x^{n-2} - \cdots - \a_1$ (since $x \ne 0$). \\
Due to $f_n(x) = x f_{n-1}(x)$, the two polynomials \f and \g have the same unique positive root $r$, that is, the graphs of the two functions \f and \g cross each other at point $r$ on the abscissa. \\
There is another intersection point between \f and \g for positive $x$ and this intersection always happens at $x = 1$. This can be easily seen from
\b
\label{at1}
f_n(1) \,\, = \,\, \a_n - \a_{n-1} - \cdots - \a_1 \,\, = \,\, f_{n-1}(1).
\e
Thus, the coordinates of the two points of intersection between \f and \g are $(r, 0)$ and $(1, y_0)$, where $y_0 = \a_n - \a_{n-1} - \cdots - \a_1$.
Next, one has to determine where point $y_0$ is with respect to the ``level" $y = \a_0 > 0$, prescribed by the free term of the Cauchy polynomial. \\
Clearly, when $y_0 = \a_n - \a_{n-1} - \cdots - \a_1 < 0$, this intersection point is in the fourth quadrant, when $y_0 = \a_n - \a_{n-1} - \cdots - \a_1 > 0$, it is in the first quadrant, and when $y_0 = \a_n - \a_{n-1} - \cdots - \a_1 = 0$, then the functions \f and \g cross only once at the abscissa at $r = 1$, i.e. $r = 1$ is a double root of \f = \g. \\
As $\a_0 > 0$, there are three possible situations: either  $0 \le y_0 \le \a_0$ (Figure 1), or $y_0 \le 0 \le \a_0$ (Figure 2), or $0 < \a_0 < y_0$ (Figure 3).
\begin{center}
\begin{tabular}{ccc}
\includegraphics[width=45mm]{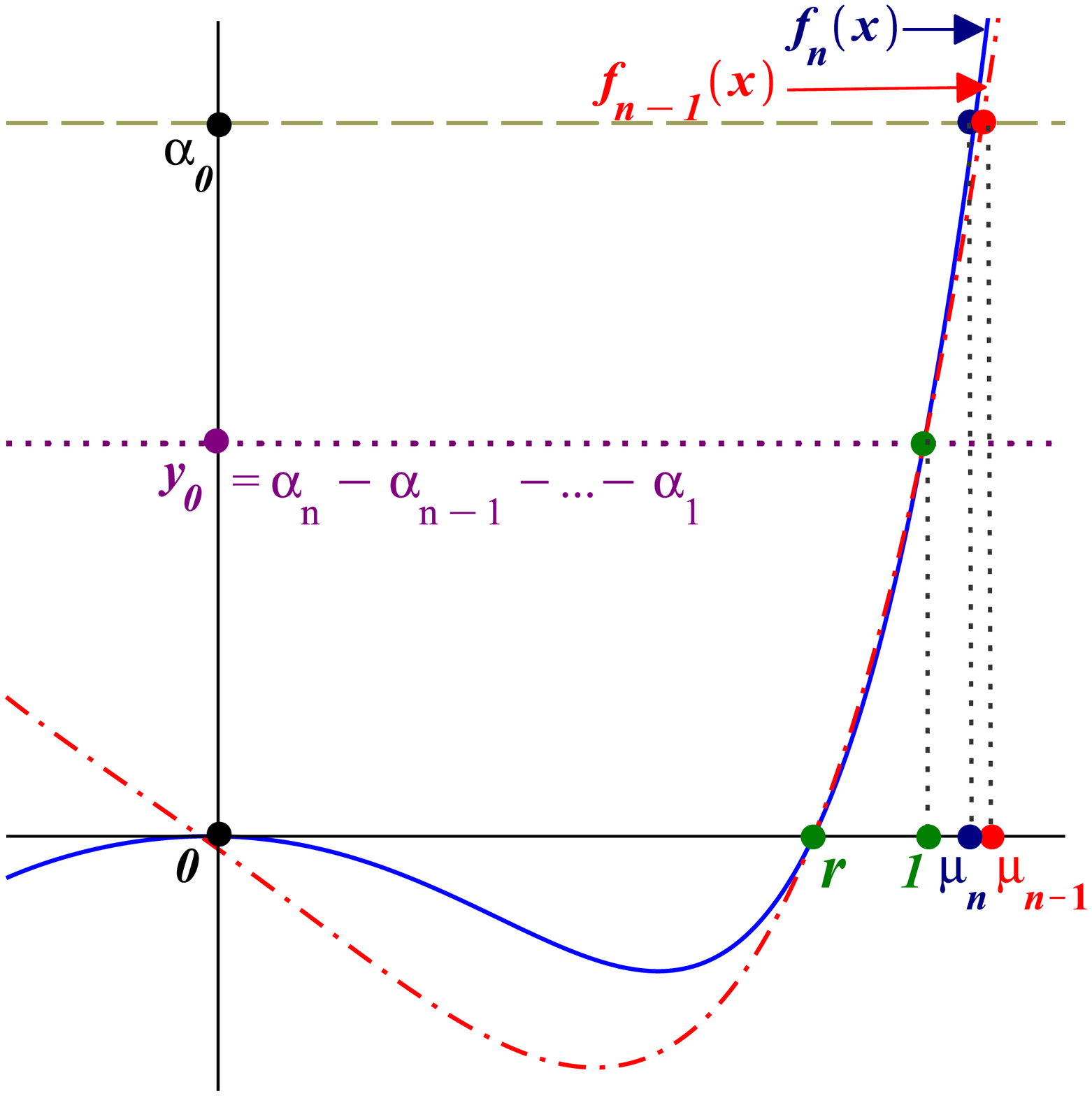} & \includegraphics[width=45mm]{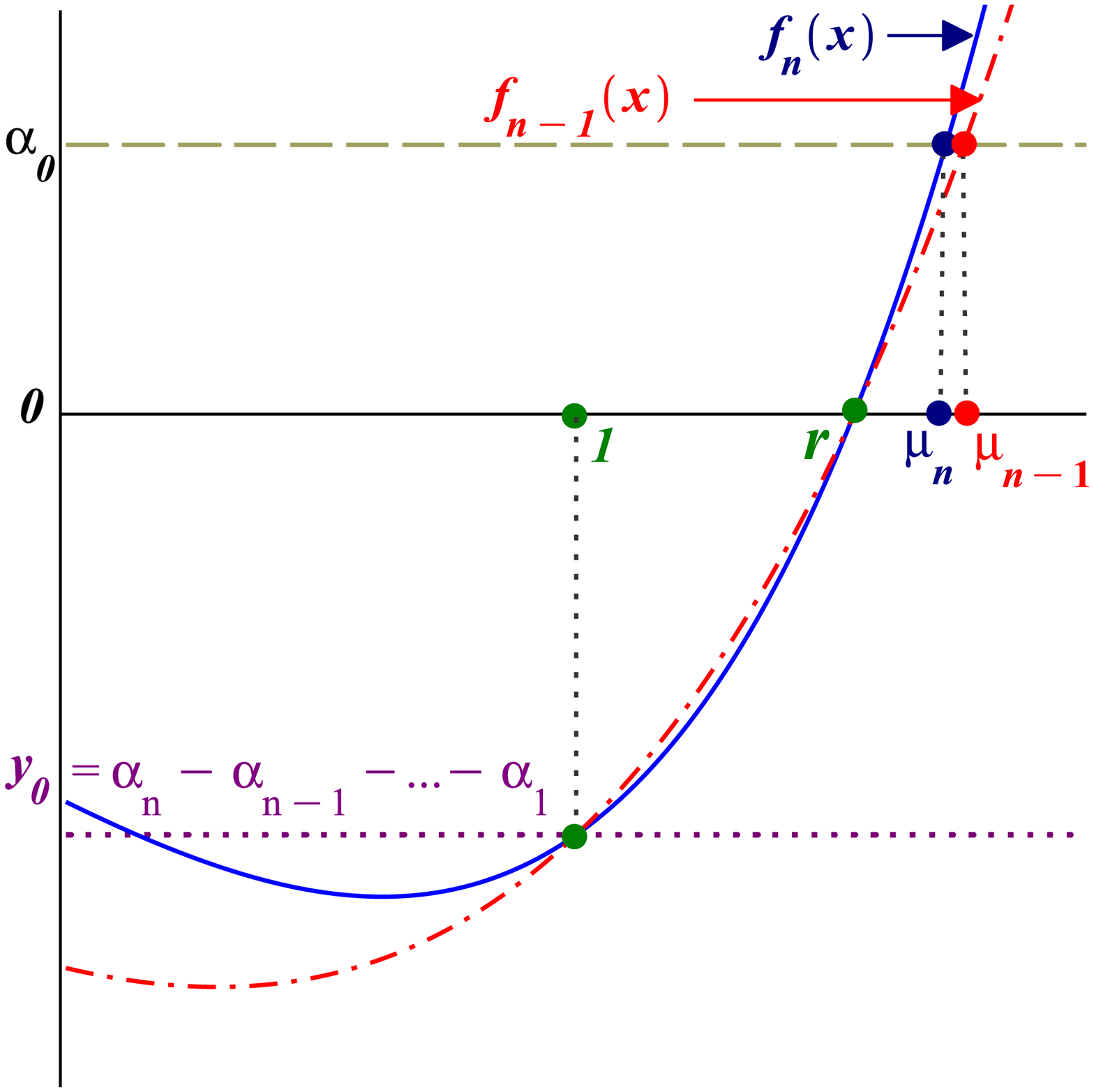} & \includegraphics[width=45mm]{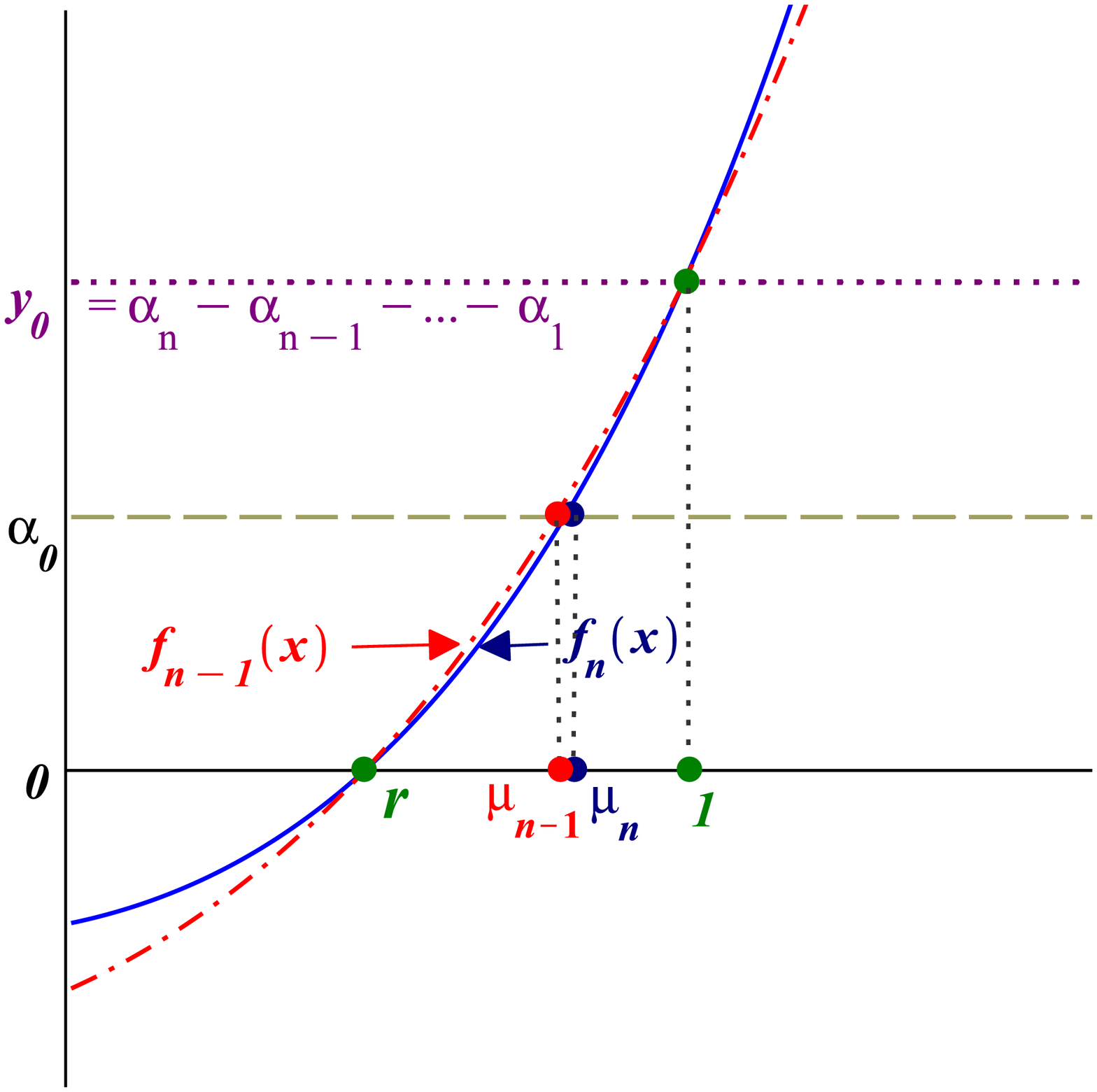} \\
{\scriptsize Figure 1} &  {\scriptsize Figure 2} & {\scriptsize Figure 3}\\
& & \\
\multicolumn{1}{c}{\begin{minipage}{11em}
\scriptsize
\begin{center}
\vskip-.3cm
\underline{$0 \le y_0 \le \a_0$}
\end{center}
When $\a_0 \ge \a_n - \a_{n-1} - \cdots - \a_1 \ge 0$, the upper bound of $\mu_n$ is $\mu_{n-1}$ and the lower bound of $\mu_n$ is 1.
\end{minipage}}
& \multicolumn{1}{c}{\begin{minipage}{11em}
\scriptsize
\begin{center}
\vskip-.3cm
\underline{$y_0 \le 0 \le \a_0$}
\end{center}
When $\a_0 \ge 0 \ge \a_n - \a_{n-1} - \cdots - \a_1$, the upper bound of $\mu_n$ is $\mu_{n-1}$ and the lower bound of $\mu_n$ is 1.
\end{minipage}}
& \multicolumn{1}{c}{\begin{minipage}{11em}
\scriptsize
\begin{center}
\vskip-.3cm
\underline{$0 < \a_0 < y_0 $}
\end{center}
When $0 < \a_0 < \a_n - \a_{n-1} - \cdots - \a_1$, the upper bound of $\mu_n$ is 1 and the lower bound of $\mu_n$ is $\mu_{n-1}$.
\end{minipage}}
\\
\end{tabular}
\end{center}
Let $\mu_n$ denote the unique positive root of \f $ = \a_0$ (note that $\mu = \mu_n$) and $\mu_{n-1}$ --- the unique positive root of \g $ = \a_0$. It is quite clear that in the first two cases (Figures 1 and 2) one has $1 < \mu_n < \mu_{n-1}$, while in the third case (Figure 3), $\mu_{n-1} < \mu_n < 1$ holds. Thus, when $\a_0 \ge \a_n - \a_{n-1} - \cdots - \a_1$, one has that $\mu_{n-1}$ is an upper bound of $\mu_n$ (or $\mu$) and 1 is a lower bound $\mu_n$ (or $\mu$). On the contrary, when $\a_0 < \a_n - \a_{n-1} - \cdots - \a_1$, one has that $\mu_{n-1}$ is the lower bound of $\mu_n$ (or $\mu$) and 1 is the upper bound $\mu_n$ (or $\mu$). In the last case, one does not need to proceed further and should just take 1 as the upper bound of the unique positive root $\mu$ of the Cauchy polynomial and, hence, as the upper bound of the roots of the general equation (\ref{eqn}). \\
If, however, the free term $\a_0$ of the Cauchy polynomial is such that $a_0 \ge \a_n - \a_{n-1} - \cdots - \a_1$, one needs to find an upper bound of $\mu_{n-1}$ and this, in turn, will serve as upper bound of $\mu_n = \mu$ and, hence, on the roots of (\ref{eqn}). This means to continue recursively by considering the equation $f_{n-1}(x) = \a_0$ and re-writing it as
\b
\label{ee}
g_{n-1}(x) = \a'_0
\e
where $g_{n-1}(x) = \a_n x^{n-1} - \a_{n-1} x^{n-2} - \cdots - \a_2 x = f_{n-1}(x) + \a_1$ and $\a'_0 = \a_0 + \a_1$, on one hand, and as
\b
\label{eee}
x g_{n-2}(x) = \a'_0,
\e
with $g_{n-2}(x) = \a_n x^{n-2} - \a_{n-1} x^{n-3} - \cdots - \a_2$, on the other hand. \\
As before, since $g_{n-1}(x) = x g_{n-2}(x)$, the polynomials $g_{n-1}(x)$ and $g_{n-2}(x)$ have the same unique positive root $r'$, that is, the graphs of the two functions $g_{n-1}(x)$ and $g_{n-2}(x)$ cross each other at point $r'$ on the abscissa. Also as before, there is another intersection point between $g_{n-1}(x)$ and $g_{n-2}(x)$ for positive $x$ and this intersection point is again $x = 1$:
\b
\label{at2}
g_{n-1}(1) \,\, = \,\, \a_n - \a_{n-1} - \cdots - \a_2 \,\, = \,\, g_{n-2}(1).
\e
Equations (\ref{ee}) and (\ref{eee}) have the same unique positive root $\mu_{n-1}$. Let $\mu_{n-2}$ denote the unique positive root of $g_{n-2}(x) = \a'_0$. This positive root is an upper bound for $\mu_{n-2}$ provided that $\a'_0 \ge \a_n - \a_{n-1} - \cdots - \a_2$. The latter is simply $a_0 \ge \a_n - \a_{n-1} - \cdots - \a_1$ and this is indeed the case, as it was presumed to hold when the recursion started. \\
The equation $g_{n-2}(x) = \a'_0$ is, in fact,
\b
\a_n x^{n-2} - \a_{n-1} x^{n-3} - \cdots - \a_3 x - \a_2 - \a_1 = \a_0.
\e
Continuing in the vein of recursively bounding the root of each of these equations with the root of an equation of degree reduced by one unit, the linear equation
\b
\a_n x - \a_{n-1} - \a_{n-2} - \cdots - \a_0 = 0,
\e
which terminates the recursion, immediately follows. The exact unique positive root of this equation is
\b
\label{moi}
\mu_1 = \frac{\a_{n-1} + \a_{n-2} + \cdots + \a_0}{\a_n}.
\e
Thus, one naturally arrives at the Cauchy bound
\b
\label{cb}
\rho = \mbox{max}\left\{ 1, \mu_1 \right\} = \mbox{max}\left\{ 1, \sum_{j = 0}^{n-1} \frac{\a_j}{\a_n} \right\}
\e
for the roots of the general equation (\ref{eqn}) --- see (8.1.10) in \cite{rs} (there are no absolute values in the above formula as the $\alpha$'s are taken as positive). \\
Next, sharper bounds than (\ref{cb}) will be found. \\
The recursion which led to (\ref{moi}) could be exited earlier. For example, if one ends at the stage of quadratic equation (preceding that of the linear equation whose root is $\mu_1$), the unique positive root $\mu_2$ of this quadratic equation will be smaller than $\mu_1$ and hence, it will provide a sharper bound. That is, a sharper bound is provided by the unique positive root of
\b
\label{q}
\a_n x^2 - \a_{n-1} x - \a_{n-2} - \cdots - \a_0 = 0,
\e
namely
\b
\mu_2 = \frac{1}{2\a_n} \left[ \a_{n-1} + \sqrt{\a^2_{n-1} + 4 \a_{n} (\a_{n-2} + \a_{n-3} + \cdots + \a_0)} \right].
\e
The recursion can be exited earlier than this --- at the stage of the cubic equation
\b
\a_n x^3 - \a_{n-1} x^2 - \a_{n-2} x - \a_{n-3} - \cdots - \a_0 = 0.
\e
Its unique positive root $\mu_3$ which is such that $\mu_3 < \mu_2 < \mu_1$, provides an even sharper bound. \\
Finally, $\mu_4$, the unique positive root of the quartic equation
\b
\a_n x^4 - \a_{n-1} x^3 - \a_{n-2} x^2 - \a_{n-3} x - \a_{n-4} - \cdots - \a_0 = 0
\e
provides the sharpest bound that can be found analytically by the recursion. \\
The Cauchy bound (\ref{cb}) can be made significantly sharper by following a different line of analysis. \\
Suppose that the number of terms with positive coefficients in the general equation (\ref{eqn}) is $p$ and that the number of terms with negative coefficients is $q$. Clearly, $p + q \le n+1$ (equality is achieved if none of the coefficients of the general equation (\ref{eqn}) is equal to zero).
For the coefficients $a_j > 0$, write $A_j$ instead, and for the coefficients $a_j < 0$, write $(-B_j)$ instead. Clearly, $A_n \equiv a_n > 0$ and all the rest of the $A$'s are non-negative. At least one of the $B$'s is positive, the rest --- non-negative (equation in which all coefficients have the same sign is no longer of interest). \\
Following the ideas, presented in \cite{1}, of splitting a polynomial into two parts and analysing the ``interaction" between these parts in order to study the roots of the polynomial, one can re-write the general equation (\ref{eqn}) as
\b
\label{forma}
A_  l x^l -  B_{n-m_1} x^{n - m_1} - B_{n-m_2} x^{n - m_2} - \cdots - B_{n-m_q} x^{n - m_q}
&& \nonumber \\
& & \hskip-8cm
= - A_{n} x^{n} - A_{n-k_1} x^{n - k_1} - \cdots - \hat{A}_  l x^l - \cdots - A_{n-k_{p-1}} x^{n - k_{p-1}},
\e
where $\{k_1 < k_2 < \ldots < k_{p-1}, m_1< m_2 < \ldots < m_q \}$ is a permutation of $\{1, 2, \ldots , n \}$,  $l$ is such that $n - m_1 < l \le n$, and the hat on $A_l$ indicates that the term $A_l x^l$ is missing from the right-hand side. At least one monomial $A_l x^l$ with $n - m_1 < l \le n$ exists --- $A_n x^n$. \\
The polynomial on left-hand side of (\ref{forma}) is a Cauchy polynomial with unique positive root $\mu$ (if the free term of the equation happens to be positive, then the resulting Cauchy polynomial will have zero as a root). This polynomial, due to having a positive coefficient in its leading term, is strictly positive for all $x > \mu$. The polynomial on the right-hand side of the equation is strictly negative for all $x > 0$ (should the free term of the original equation happen to be negative, the polynomial on the right-hand side will have 0 as a root). For all $x > \mu$, the two sides have opposite sign and therefore, there can be no roots of the equation for $x > \mu$, i.e. the root bound for the general equation (\ref{eqn}) is $\mu$ --- the unique positive root of the Cauchy polynomial extracted from the given polynomial. Of course, different choices of $l$ in $A_l x^l$ will lead to different Cauchy polynomials  with different roots on the left-hand side of (\ref{forma}). It will be the unique positive root of the particular Cauchy polynomial appearing on the left-hand side of (\ref{forma}) that will provide a bound for the roots of the general equation (\ref{eqn}). \\
Suppose now that on the left-hand side of (\ref{forma}) the term $-B_{l-1}x^{l-1} - B_{l-2}x^{l-2} - \cdots - B_{m-n_1+1}x^{n-m_1+1}$ with $B_{l-1} = B_{l-2} = \ldots = B_{m-n_1+1} = 0$ has been added to the Cauchy polynomial. Following the analysis that lead to the derivation of (\ref{moi}) and (\ref{cb}), it can be seen that, instead of the bound given by the Cauchy formula (\ref{cb}), a significantly lower bound can be used
\b
\label{new}
\rho'_l = \mbox{max}\left\{ 1, \sum_{j = 1}^{q} \frac{B_{n-m_j}}{A_l} \right\}.
\e
The summation is no longer over all coefficients of the equation (as in the Cauchy formula), but only over the absolute values of the negative ones (in units of $A_l$). There will be different bounds $\rho'_l$ for different choices of $l$. \\
Clearly, the sharpest of these bounds will be the one with the biggest denominator: $A_{\mbox{\tiny max}} = $ max$\{A_l \, | \, n \ge l > n-m_1 \}$:
\b
\rho' =  \mbox{max}\left\{ 1, \sum_{j = 1}^{q} \frac{B_{n-m_j}}{A_{\mbox{\tiny max}}} \right\}.
\e
It should also be noted that, instead of introducing the ``ghost" term $B_{l-1}x^{l-1} - B_{l-2}x^{l-2} - \cdots - B_{m-n_1+1}$ with $B_{l-1} = B_{l-2} = \ldots = B_{m-n_1+1} = 0$, a sharper bound can be obtained if one terminates the recursion when the equation
\b
A_l x^{l-n-m_1} - B_{n-m_1} - B_{n-m_2} - \cdots - B_{n-m_q} = 0
\e
is reached. The obtained in such way bound, for different values of $l$, is
\b
\label{new2}
\rho''_l = \mbox{max}\left\{ 1, \left( \sum_{j = 1}^{q} \frac{B_{n-m_j}}{A_l} \right)^{\frac{1}{l-n-m_1}} \right\}, \,\,\,\, n-m_1 < l \le n.
\e
From here, one can find $\rho''$ as the minimum of the above. \\
The above arguments prove the following theorem
\3
\begin{center}
\begin{minipage}{13cm}
{\it An upper bound of the roots of the general polynomial $a_n x^n + a_{n-1} x^{n-1} + \cdots + a_0$ with $a_n > 0$ is the smallest of the unique positive roots of all Cauchy polynomials that can be extracted form this polynomial with preservation of all terms with negative coefficients.}
\end{minipage}
\end{center}
\3
\n
To determine the lower bound of the unique positive root of the Cauchy polynomial, re-write equation (\ref{stump}) for positive $x$ as
\b
x^n \left( \a_n  - \a_{n-1} \frac{1}{x} - \cdots - \a_0 \frac{1}{x^n} \right) = 0.
\e
In variable $y = 1/x$, this equation becomes:
\b
\label{pp}
\a_0 y^n + \a_1 y^{n-1} + \cdots + \a_{n-1} y - \a_n = 0.
\e
As there are $n$ Cauchy polynomials $\a_{n-t} y^t - \a_n$ ($t = 1,2, \ldots, n$) that can be extracted from $\a_0 y^n + \a_1 y^{n-1} + \cdots + \a_{n-1} y - \a_n$, equation (\ref{pp}) can be written in $n$ equivalent ways:
\b
\label{extra}
- \a_n + \a_{n-t} y^t = - \a_0 y^n - \a_1 y^{n-1} - \cdots - \hat{\a}_{n-t} y^t - \cdots - \a_{n-1} y.
\e
For each of these, the unique positive root $(\a_{n} / \a_{n-t})^{1/t}$ of the Cauchy polynomial $\a_{n-t} y^t - \a_n$ provides an upper bound on the roots of equation (\ref{pp}). The full set of the obtained in this way bounds on the roots of this equation is
\b
\label{sq}
\left\{ \frac{\a_n}{\a_{n-1}}, \left( \frac{\a_n}{\a_{n-2}} \right)^{\frac{1}{2}}, \ldots, \left( \frac{\a_n}{\a_0} \right)^{\frac{1}{n}} \right\}.
\e
Thus, the sharpest upper bound on the roots of equation (\ref{pp}) is, obviously, the smallest of all these numbers. This is a Lagrange type of bound (the Lagrange bound is the sum of the two largest values of this sequence). \\
The sharpest lower bound of the unique positive root of the Cauchy equation (\ref{stump}) will  thus be the largest of the reciprocals of the numbers in the sequence (\ref{sq}).

\end{document}